\documentclass[12pt]{article}
\usepackage[dvips]{epsfig}
\usepackage[dvips]{graphics}
\usepackage[dvips]{color}
\usepackage{amsfonts}
\usepackage{amsmath}
\usepackage{amssymb}
\usepackage{amsthm}
\usepackage{graphics}

\newtheorem{thm}{Theorem}
 \newtheorem{cor}[thm]{Corollary}
 
 \newtheorem{lem}[thm]{Lemma}
 \newtheorem{prop}[thm]{Proposition}

\newtheorem{defn}[thm]{Definition}
 \newtheorem{rem}{Remark}

\begin{document}

\begin{center}
\textbf{\Large{ On The Slope of Hyperelliptic Lefschetz Fibrations
and The Number of Separating Vanishing Cycles \\
 }}
 \bigskip
Yusuf Z Gurtas
\end{center}
\begin{abstract} In this article we find an upper bound for the slope of genus $g$
hyperelliptic Lefschetz fibrations, which is sharp when $g=2$, and
demonstrate the strong connection, in general,    between the slope
of hyperelliptic genus $g$ Lefschetz fibrations  and the number of
separating vanishing cycles. Specifically, we show that the slope is
greater than $4-\frac{4}{g}$ if and only if the fibration contains
separating vanishing cycles. We also improve the existing bound on
$\frac{s}{n}$, the ratio of number of separating vanishing cycles to
the number of non-separating vanishing cycles, for hyperelliptic
Lefschetz fibrations of genus $g\geq2$. In particular we show that
$s\leq n$ for such fibrations when $g\geq 6$.
\end{abstract}

\author{Yusuf Z Gurtas}

\section{Introduction} 
\label{intro}

Let $X \rightarrow S^{2}$ be a genus $g$   Lefschetz fibration. (The
reader is referred to \cite{GS} for a thorough review of Lefschetz
fibrations.)

 It's known that the $4-$ manifold $X$ carries an almost
complex structure; therefore it makes sense to define its
holomorphic Euler characteristic and first Chern class. Let
\[\chi _{h}=\frac{1}{4}\left( \sigma +\chi \right)
 \ \ \textup{and} \ \ c_{1}^{2}=2\chi +3\sigma,\] where $\chi$ is the Euler
characteristic and $\sigma$ is the signature of the $4-$ manifold
$X$. The \emph{slope} $\lambda _{f}$ of $X$ is defined as $\lambda
_{f}:=K_{f}^{2}/\chi _{f}$ where $K_{f}^{2}:= c_1^2+8(g-1) $ and
$\chi _{f}:=\chi _{h}+g-1.$  It's known that
\[\lambda _{f}\geq 4-\frac{4}{g}\] for a genus $g$ Lefschetz
fibration   and this bound is sharp. For example, all of the known
hyperelliptic Lefschetz fibrations over $S^{2}$ with no separating
vanishing cycles satisfy $\lambda _{f}= 4-\frac{4}{g}.$ We will
write $\lambda$ for simplicity from now on and all the Lefschetz
fibrations discussed in this article will be hyperelliptic.

 The connection between $\lambda$ and the number of
separating vanishing cycles of a Lefschetz fibration seems to be
unaccounted for in the literature. Let $s$ be the number of
separating vanishing cycles and $n$ be the number of those that are
non-separating. In this article we will prove:
\begin{thm} \label{main-thm} A genus $g$ hyperelliptic
Lefschetz fibration $X \rightarrow S^{2}$ satisfies $\lambda
> 4-\frac{4}{g}$ if and only if $s\neq0$, i.e., it contains separating vanishing
cycles.
\end{thm}
Recall that  a Lefschetz fibration can not contain only separating
vanishing cycles. Therefore the theorem should be interpreted as a
fibration containing a mixture of separating and non-separating
vanishing cycles.

 An interesting question that arises at this
point is the proportion of the number of separating cycles within a
fibration, in particular its ratio to the number of non-separating
vanishing cycles, \  $ \frac{s}{n}$ .  We   do not find  any
estimates in the literature on this ratio except for
\begin{eqnarray} \label{Stipsicz-bound}
\frac{s}{n}&\leq &5
\end{eqnarray}
due to A.Stipsicz, \cite{S}. Since we have $n>0$ in a given
Lefschetz fibration, this ratio is always defined.
\begin{defn}
\begin{eqnarray*}
\rho \left( g\right) =\max \left\{  \right.  r= \frac{s}{n} &|&\!
\exists \, \textup{ a
}\Sigma_{g}-\textup{hyperelliptic Lefschetz fibration } X \rightarrow S^{2} \textup{ with } s\\
&& \!\! \left.  \textup{ separating and } n \ \textup{non-separating
vanishing cycles} \right\}
\end{eqnarray*}
\end{defn}There isn't enough evidence to justify that the  bound
(\ref{Stipsicz-bound}) could actually be sharp. On the contrary, all
of  the known examples suggest that   $\rho\left( g\right)$  may not
be too high.

In this article we will improve the bound on $\rho$  for
hyperelliptic Lefschetz fibrations and  show that:

\begin{thm} \label{ratio-estimate}   For an hyperelliptic    Lefschetz fibration of genus $g\geq2$ we have
\[ \rho \left( g\right) \leq \frac{3g+2}{4\left( g-1\right) }.\]
\end{thm}

The last result is about signature of hyperelliptic Lefschetz
fibrations.

Even though there is an explicit formula that gives the signature in
terms of separating and non-separating vanishing cycles
 for genus  $g$ hyperelliptic Lefschetz fibrations,
 it is desirable to have a formula that relates the signature to the total number
of vanishing cycles, perhaps by a scalar multiplication.
\begin{thm} \label{signature-in-terms-of-lambda} For a genus  $g$ hyperelliptic Lefschetz fibration we have
 \[\sigma =k\left( n+s\right),\]
 where $\displaystyle k= -\allowbreak \frac{\lambda-8}{\lambda-12}.$
\end{thm}

In the next section we will prove Theorem \ref{main-thm} and Theorem
\ref{signature-in-terms-of-lambda} and show some of their
applications for genus $2$. The following section will summarize
similar results for genus $3$. The case of low genus is handled
separately  because there is only one type of separating vanishing
cycle when $g<4$ and due to that reason general formulas don't
always give rise to results that are as sharp as  could be when
restricted to low genus. It is also intended to give the reader an
easy preparation for the general case which will be addressed in the
last section along with the proof of Theorem \ref{ratio-estimate}.\\

We prove all the results for heyperelliptic Lefschetz fibrations but
some of them generalize to non-heyperelliptic case as well. Please
see Remark \ref{generalizations-rem} for results that generalize to
non-heyperelliptic Lefschetz fibrations. Even though we found out
 that there are  shorter proofs for some of the results, we chose to
 leave them in the original format they were written in. We pointed
 out to those shorter proofs in Remark
 \ref{final-remarks-and-alternate-proofs}.
 We don't claim originality on most of the results but Theorem
\ref{ratio-estimate} has not appeared anywhere else to the best of
our knowledge.

\section{Genus $2$ }

The signature of a genus $g$ hyperelliptic Lefschetz fibration $X
\rightarrow S^{2}$ is given by
\[-\frac{g+1}{2g+1}n+\sum_{h=1}^{\left[ g/2\right] }\frac{4h\left( g-h\right)
s_{h}}{2g+1}-s.\] Let \[x=\sum_{h=1}^{\left[ g/2\right] }h\left(
g-h\right) s_{h},\] where $\displaystyle s=\sum_{h=1}^{\left[
g/2\right] }s_{h}.$ The other invariants of
$X$ that will be used throughout the article are:\\
  Euler characteristic
\[\chi =n+s-4\left( g-1\right),\]
holomorphic Euler characteristic
\begin{eqnarray} \label{holom-in-terms-of-x}
\nonumber \chi_{h}=\frac{1}{4}\left( \chi +\sigma \right)
&=&\frac{1}{4}\left( n+s-4\left( g-1\right)
-\frac{g+1}{2g+1}n+\frac{4x}{2g+1}-s\right)\\
&=& \frac{ng+4x }{4\left( 2g+1\right)} -(g-1),
\end{eqnarray}
and square of the first Chern class $c_1^2$
\begin{eqnarray*}
 c_{1}^{2}=2\chi +3\sigma &=&2\left( n+s-4\left( g-1\right) \right) +3\left( -\frac{g+1}{2g+1}n+\frac{4x}{2g+1}-s\right)\\
&=& 2n-s-8\left( g-1\right) -3\frac{g+1}{2g+1}n+\frac{12x}{2g+1},
\end{eqnarray*}
where $s$ is the number of separating vanishing cycles and $n$ is
the number of non-separating vanishing cycles.
\begin{lem} \label{sg-less-than-2x} $sg\leq 2x$ for $g\geq 2.$
\end{lem}
\begin{proof} It's not difficult to see that $s(g-1)\leq x$ by definition of $x$ and $s$.
Therefore
\begin{eqnarray*}  s\leq \frac{x}{g-1} \ \ \ \textup{and} \ \ \ sg\leq \frac{gx}{g-1}.
\end{eqnarray*}
 The proof follows from the fact that $\frac{g}{g-1}\leq 2$ for $g\geq 2.$
\end{proof}
\begin{proof}[Proof of Theorem  \ref{main-thm}]
The slope $\lambda $ of the fibration is given as
\begin{eqnarray} \nonumber
\lambda &=&\frac{c_{1}^{2}+8\left( g-1\right) }{\chi _{h}+\left(
g-1\right) } = \frac{2n-s  -3\frac{g+1}{2g+1}n+\frac{12x}{2g+1}
}{\frac{ng+4x }{4\left( 2g+1\right)}}\\
\label{slope-formula} &=&4\frac{n\left( g-1\right) -s\left(
2g+1\right) +12x}{ng+4x}\,.
\end{eqnarray}
Assume $s\neq 0$. Then $x\neq0$ and we have
\begin{eqnarray*}
\lambda - (4-4/g)&=&4\frac{n\left( g-1\right) -s\left( 2g+1\right)
+12x}{ng+4x}-4+4/g\\
&=&4\frac{-2sg^{2}-sg+8gx+4x}{\left( ng+4x\right) g} =4\frac{\left(
2g+1\right) \left( 4x-sg\right) }{\left( ng+4x\right) g}> 0,
\end{eqnarray*}
because $4x> sg$ by Lemma \ref{sg-less-than-2x}. and all other
factors are positive. Therefore $\lambda - (4-4/g)=0$ if and only if
$4x=sg$; i.e., if and only if $s=0$.
\end{proof}

\begin{cor}  \label{genus2-slope} For a genus $2$ Lefschetz fibration we have
\begin{equation} \label{slope-in-n-and-s} \lambda =2\,\frac{n+7s}{n+2s}=2\,\frac{1+7r}{1+2r}.
\end{equation}
\end{cor}

\begin{proof}
 From (\ref{slope-formula})  we have
\[\lambda= 4\frac{n\left( g-1\right) -s\left( 2g+1\right)
+12x}{ng+4x}\ . \] Setting $g=2$ and realizing that for a genus $2$
Lefschetz fibration $x=s$ we obtain
\[\allowbreak 4\frac{n\left( 2-1\right) -s\left( 2\cdot 2+1\right)
+12s}{2n+4s}=2\frac{n+7s}{n+2s}.\] Dividing through by $n$ gives
\[2\,\frac{1+7r}{1+2r}. \]
\end{proof}

\begin{prop}[(Corollary 10, \cite{Oz} )] \ \label{prop-burak} For a genus $2$ Lefschetz fibration we have
\[c_{1}^{2}\leq 6\chi _{h}-3.\]
\end{prop}

\begin{proof}
We will use the bound

\begin{equation} \label{burak-sign-bound} \sigma \leq n-s-4
\end{equation} for hyperelliptic Lefschetz fibrations given by Corollary 9, \cite{Oz}.
First, we write $\chi $ in terms of $\chi_{h}$:
\begin{eqnarray*}
\chi _{h}&=&\frac{1}{4}\left( \sigma +\chi \right)=\frac{1}{4}\left(
-\frac{3}{5}n-\frac{1}{5}s+n+s-4\right)
=\frac{1}{10}n+\frac{1}{5}s-1, \\\ \ \textup{therefore} \\ \chi
&=&n+s-4=10\left( \frac{1}{10}n+\frac{1}{5}s-1\right)
+6-s=10\chi _{h}+6-s.\\
\textup{Then, since}\ \ \ \sigma &\leq&
n-s-4=n+s-4-2s=\chi -2s, \ \ \textup{we have} \\
\chi _{h}&=&\frac{1}{4}\left( \sigma +\chi \right)  \leq
\frac{1}{4}\left( \chi -2s+\chi \right)=\frac{1}{2}\left( \chi
-s\right)\\\chi_h&\leq & \frac{1}{2}\left( 10\chi _{h}+6-s-s\right)
=5\chi _{h}+3-s,
\end{eqnarray*}
which can be written as
\begin{equation} \label{s-and-chi_h} s\leq 4\chi _{h}+3.
\end{equation}
Finally, we have
\begin{eqnarray} \label{Noether}
c_{1}^{2}=12\chi _{h}-\chi& =&12\chi _{h}-\left( 10\chi
_{h}+6-s\right) =2\chi _{h}-6+s\\ \nonumber &\leq &2\chi
_{h}-6+4\chi _{h}+3=6\chi _{h}-3.
\end{eqnarray}
\end{proof}

\begin{rem} \label{holom-nonneg} \textup{ Solving the   inequality (\ref{s-and-chi_h}) for $\chi _{h}$ we get \
$\displaystyle \frac{1}{4}\left( s-3\right) \leq \chi _{h}$. This
means that for genus $2$ Lefschetz fibrations we have $\chi_{h}\geq
0.$ }
\end{rem}

\begin{rem} \textup{ For  genus $2$ Lefschetz fibrations we have $c_{1}^{2}=2\chi
_{h}+s-6$ by (\ref{Noether}). Therefore all  genus $2$ fibrations
with no separating vanishing cycles are necessarily on the Noether
line. The manifold lands above Noether line if and only if it
contains separating vanishing cycles.}
\end{rem}
\begin{cor}  \label{slope-estimate}For a genus  $2$ Lefschetz fibration we have \[\lambda \leq 6-\frac{1}{\chi
_{h}+1}.\]
\end{cor}

\begin{proof}
Using Proposition  \ref{prop-burak} we can write \[c_{1}^{2}+8\leq
6\chi _{h}-3+8=6\left( \chi _{h}+1\right) -1.\]Dividing through by
$\chi _{h}+1$ we obtain \[\lambda =\frac{c_{1}^{2}+8}{\chi
_{h}+1}\leq 6-\frac{1}{\chi _{h}+1}.\] Note that $\chi _{h}+1>0$ by
Remark  \ref{holom-nonneg}.
\end{proof}

\begin{cor}  \label{genus-2-rho} $ \rho\left( 2\right) \leq 2. $
\end{cor}

\begin{proof}
Using Corollary \ref{genus2-slope} and Corollary
\ref{slope-estimate} we can write
\[\lambda=2\,\frac{1+7r}{1+2r}\leq 6-\frac{1}{\chi _{h}+1}\leq6.\]
for any genus  $2$ Lefschetz fibration. Solving it for $r$ gives
$r\leq 2$.
\end{proof}

\begin{cor} \label{system-of-equality-inequality} The number of separating and non-separating vanishing
cycles $s$ and $n$, respectively, in a  genus  $2$ Lefschetz
fibration satisfy
\begin{eqnarray*}2s+n&=&10k\\ 2n-s&\geq&5
\end{eqnarray*}
for some $k\in \mathbb{Z}^+.$
\end{cor}
\begin{proof}
\begin{eqnarray}
\chi_h=\frac{1}{4}\left(\sigma +\chi\right) &=&\frac{1}{4}\left(-\frac{3}{5}n-\frac{1}{5}s+n+s-4\right)\\
&=&\frac{1}{10}\left( n+2s\right) -1 \label{holom-in-n-and-s}
\end{eqnarray}
Therefore $n+2s=10\left(\chi_h+1\right)$ and $\chi _{h}+1>0$ by
Remark \ref{holom-nonneg}. This proves the equality. For the
inequality we will use Corollary \ref{genus2-slope} and Corollary
\ref{slope-estimate} :
\begin{eqnarray*}
\ 2\frac{n+7s}{n+2s}\leq 6-\frac{1}{\chi
_{h}+1}&=&6-\frac{1}{\frac{1}{10}\left( n+2s\right)
-1+1}=6-\frac{10}{n+2s}.
\end{eqnarray*}
Solving \[  2\frac{n+7s}{n+2s} \leq 6-\frac{10}{n+2s} \] for $s$ we
obtain
 $s\leq 2n-5$ as claimed.
\end{proof}
\indent It would be an interesting question to ask if this
inequality
 is  sharp.
 \begin{prop} If the equations
 \begin{eqnarray*}2s+n&=&10k\\ 2n-s&=&5
\end{eqnarray*}
are satisfied for a genus  $2$ Lefschetz fibration  then
 \begin{eqnarray} \label{slope-using-m} \frac{s}{n}=\frac{4m+3}{2m+4}\end{eqnarray}
for $m\geq 0$.
 \end{prop}

 \begin{proof}
 Solving the given system of equations we
obtain \[n=2+2k,s=-1+4k\] $k\in \mathbb{Z}^+.$ Therefore
\[\frac{s}{n}=\frac{4k-1}{2k+2}.\] Now, let $m=k-1\geq0 $.
 \end{proof}

 First few values this sequence can take on are
\[\frac{s}{n}=\frac{3}{4},\frac{7}{6},\frac{11}{8},\frac{3}{2},\frac{19}{12},\frac{23}{14}.\]
Xiao constructed examples realizing the values
$\frac{3}{4},\frac{7}{6}$ and $\frac{19}{12}$, \cite{X}.
\begin{rem} \textup{ With $\frac{s}{n}=\frac{4m+3}{2m+4}$ the
slope becomes:
\[\lambda=2\frac{1+7\frac{4m+3}{2m+4}}{1+2\frac{4m+3}{2m+4}}=\frac{6m+5}{m+1}=6-\frac{1}{m+1}\ .\]
Invoking Corollary \ref{slope-estimate} we get
\[6-\frac{1}{m+1}\leq 6-\frac{1}{\chi _{h}+1}, \ \ \textup{i.e.,}\ \ 0\leq m\leq\chi
_{h}.\] }
\end{rem}

It's interesting to note that this bound is sharp for the  examples
that we know satisfy the equation $2n-s=5$, i.e., $m=\chi_h$.
Therefore we might conjecture that this is a characterizing feature
for genus $2$ fibrations satisfying $2n-s=5$. Indeed that is the
case:
\begin{prop}\label{bound-on-slope-is-sharp-prop}   For a genus $2$ Lefschetz
fibration we have
\[ 2n-s=5 \ \ \ \textup{if and only if } \ \ \ \lambda=6-\frac{1}{\chi _{h}+1}. \]
\end{prop}

\begin{proof}
 Assume  $2n-s=5$. Substitute $\chi _{h}=\frac{1}{10}\left( n+2s\right)
 -1$ into
 \begin{equation} \label{ratio-in-terms-of-holo}
 \frac{4\chi _{h}+3}{2\chi _{h}+4}
 \end{equation}
  and use $2n-s=5$ for both the numerator and denominator
 to see that it's equal to $\frac{s}{n}$. Then substitute
 (\ref{ratio-in-terms-of-holo}) in place of $\frac{s}{n}$ in (\ref{slope-in-n-and-s})  to obtain the desired equality.
 Conversely, assume that the bound on $\lambda$ is sharp. Substitute   $\chi _{h}=\frac{1}{10}\left( n+2s\right)
 -1$ into the bound and set it equal to (\ref{slope-in-n-and-s}).
 Solving that equality for $s$ will result in  $s=2n-5$.
\end{proof}

\begin{rem} \label{invariants-for-region-III}\textup{
We calculate the invariants of a genus $2$ Lefschetz fibration with
$2n-s=5$  as :} \vspace{-.1in}
\begin{eqnarray*}
\sigma &=&-n+1=-\frac{1}{2}\left( s+3\right) \\
\chi &=&\allowbreak 3n-9=\frac{3}{2}\left( s-1\right) \\
\chi _{h}&=&\allowbreak \frac{1}{2}n-2=\frac{1}{4}\left( s-3\right)
\\
c_{1}^{2}&=&\allowbreak 3n-15=\frac{3}{2}\left( s-5\right)
\end{eqnarray*}
 \end{rem}
 \begin{rem}\textup{ The bound (\ref{burak-sign-bound})   on signature is sharp
 and  realized by genus $2$ Lefschetz fibrations satisfying $2n-s=5$. Simply write
 $2n-s=5$ as $ n-s-4 = -n+1=\sigma.$ }
 \end{rem}
 \begin{rem} \label{slope-in-n-and-s-only} \textup{ Thanks to the computations in
 Remark \ref{invariants-for-region-III} we can express the slope $\lambda$ in terms of $n$
 and $s$ only as
\begin{eqnarray*}\lambda =\frac{3n-15+8}{\allowbreak
\frac{1}{2}n-2+1}=\allowbreak
 2\frac{3n-7}{n-2} \ \ \ \textup{and} \ \ \
 \lambda =\frac{\frac{3}{2}\left( s-5\right) +8}{\frac{1}{4}\left( s-3\right)
 +1}=2\frac{1+3s}{1+s},
 \end{eqnarray*}
 respectively, for fibrations satisfying $2n-s=5$. }
 \end{rem}
  Combining the results on the slope of genus  $2$ Lefschetz fibrations
   so far with Proposition \ref{slope-estimate-general-prop}
  and Propositioin \ref{slope-double-estimate-prop} we can prove:
\begin{cor} \label{slope-bound-in-n-and-s-only}
For a genus  $2$ Lefschetz fibration  with $n$ non-separating and
$s$ separating vanishing cycles we have
\begin{eqnarray*} \lambda =2\frac{n+7s}{n+2s}\leq
2\frac{1+3s}{1+s}\leq 2\frac{6s+3n-5}{2s+n}\leq 2\frac{3n-7}{n-2}
\leq 10\frac{s+n-2}{2s+n} \leq 2\frac{5n-s-12}{n-2}.
\end{eqnarray*}
\end{cor}
\begin{proof}
All but the fourth inequality are equivalent to  $2n-s\geq 5$, which
is true by Corollary  \ref{system-of-equality-inequality}. The
fourth inequality turns out to be $0\leq 2\left( n-4\right) \left(
2n-s-5\right)$ but this is also true thanks to Corollary
\ref{system-of-equality-inequality} and Remark
\ref{multiple-of-five}.
All five inequalities become equality when $2n-s=5$.
 \end{proof}

Now, we will prove Theorem \ref{signature-in-terms-of-lambda}.

\begin{proof}[Proof of Theorem \ref{signature-in-terms-of-lambda}] From (\ref{slope-formula}) we have
\[\lambda= 4\frac{n\left( g-1\right) -s\left( 2g+1\right)
+12x}{ng+4x}\ . \] Cross multiplication   gives
\[4n\left( g-1\right) -4s\left( 2g+1\right) +48x = \lambda ng+4x\lambda.\]
Solving this for $x$ results in \[x = \frac{\lambda ng+4s\left(
2g+1\right) -4n\left( g-1\right) }{4\left( 12-\lambda\right) }.\]
 We will substitute this into the signature formula to obtain
the result:
\begin{eqnarray*}
\sigma &=&-\frac{g+1}{2g+1}n+\sum_{h=1}^{\left[ g/2\right]
}\frac{4h\left( g-h\right) s_{h}}{2g+1}-s \\
&=&-\frac{g+1}{2g+1}n+\frac{4x}{2g+1}-s\\
&=& -\frac{g+1}{2g+1}n+\frac{\lambda ng+4s\left( 2g+1\right)
-4n\left(
g-1\right) }{\left( 2g+1\right) \left( 12-\lambda\right) }-s \\
&=&\frac{-\left( g+1\right) n\left( 12-\lambda\right) +\lambda
ng+4s\left( 2g+1\right) -4n\left( g-1\right) -s\left( 2g+1\right)
\left(
12-\lambda\right) }{\left( 2g+1\right) \left( 12-\lambda\right) }\\
&=&\frac{\left( 2g+1\right) \left( \lambda-8\right) \left(
n+s\right)
}{\left( 2g+1\right) \left( 12-\lambda\right) }\\
&=&-\frac{8-\lambda}{12-\lambda}\left( n+s\right).
\end{eqnarray*}
\end{proof}
\begin{rem} \label{signature-negative}  \textup{ We have $c_{1}^{2} <  8\chi _{h}$ for  genus  $2$ Lefschetz
fibrations by Proposition \ref{prop-burak}. Therefore $\lambda <8$
and the signature is always
  negative for those fibrations by Theorem
  \ref{signature-in-terms-of-lambda}. }
\end{rem}

\begin{cor} \label{cor-signature-holom-bound} For a genus  $2$ Lefschetz fibration we have
\begin{eqnarray} \label{signature-holom-bound} \sigma \leq -2\chi _{h}-3
 \ \ \ \textup{and} \ \ \ \sigma \leq -\frac{1}{3}\chi -2.
 \end{eqnarray}
\end{cor}
\begin{proof}
We have $\displaystyle \lambda \leq 6-\frac{1}{\chi
 _{h}+1}$ by Corollary \ref{slope-estimate} and $\displaystyle -\frac{8-\lambda}{12-\lambda}$ is
 an increasing function of $\lambda$. Therefore, substituting
 $\displaystyle 6-\frac{1}{\chi _{h}+1}$ in place of $\lambda$
 in Theorem \ref{signature-in-terms-of-lambda} gives \[\sigma \leq
 -\frac{8-\left( 6-\frac{1}{\chi _{h}+1}\right) }{12-\left( 6-\frac{1}{\chi _{h}+1}\right) }\left( n+s\right)  \leq -\frac{2\chi _{h}+3}{6\chi
_{h}+7}\left( n+s\right)=-\frac{2\chi _{h}+3}{6\chi _{h}+7}\left(
\chi+4\right).\]  Now, substitute $\chi=4\chi_h-\sigma$ and cross
multiply to get \[\sigma \left( 6\chi _{h}+7\right) \leq -\left(
2\chi _{h}+3\right) \left( 4\chi _{h}-\sigma +4\right) \] using
$\chi_h\geq0$. Solving this for $\sigma$ gives   the first
inequality. In order to obtain the second inequality simply
substitute $\chi_h=\frac{1}{4}\left( \chi +\sigma \right) $ into the
first one and solve for $\sigma$. Note that both inequalities are
sharp for  genus $2$ fibrations with $2n-s=5$ and they can also be
obtained using Remark \ref{invariants-for-region-III} in that case.
 \end{proof}
 \begin{rem} \label{multiple-of-five} \textup{ We proved in Corollary \ref{system-of-equality-inequality} that $2n-s\geq 5$ for
 genus $2$ Lefschetz fibrations.
 In fact $2n-s$ is divisible by $5$:
 \begin{eqnarray*} 2n-s&=&2\left( n+2s\right) -5s=20\left( \chi_h+1\right) -5s=
 5\left( \allowbreak 4\chi_h+4-s\right)  \\
 &=&5\left( \allowbreak \sigma +\chi +4-s\right)=5\left( \allowbreak n+\sigma \right).
 \end{eqnarray*}
 (One can also use the local signature formula $\sigma =-\frac{3}{5}n-\frac{1}{5}s$ in order to see
 that, \cite{Ma})
 Let $t=n+\sigma.$ It's clear that $t\in\mathbb{Z}^{+}.$ Solving the equations
 \begin{eqnarray*}2s+n&=&10k\\ 2n-s&=&5t
\end{eqnarray*}
for $n$ and $s$ we get $n=2t+2k,s=4k-t$. In particular $n\geq 4$
because $t,k\in\mathbb{Z}^{+}$. Substituting these values of $n$ and
$s$ in (\ref{slope-in-n-and-s}) we obtain
\[\lambda=2\frac{1+7\frac{4k-t}{2t+2k}}{1+2\frac{4k-t}{2t+2k}}=\allowbreak 6-\frac{t}{k}=\allowbreak 6-\frac{t}{\chi _{h}+1}
\leq \allowbreak 6-\frac{1}{\chi _{h}+1}\]
 as we proved in Corollary \ref{slope-estimate}. }
\end{rem}
 \begin{cor} \label{proof-2-for-genus-2-bound} For genus $2$ Lefschetz fibrations we have
\begin{eqnarray*}\frac{s}{n}\leq \frac{4\chi _{h}+3}{2\chi _{h}+4}.
\end{eqnarray*}
\end{cor}
\begin{proof} Using (\ref{holom-in-n-and-s}) and $s\leq 2n-5$ we have
\begin{eqnarray*}\chi _{h}=\frac{1}{10}\left( 2s+n\right) -1\leq \frac{1}{10}\left( 2\left( 2n-5\right) +n\right) -1=\frac{1}{2}n-2.
\end{eqnarray*}
Thus $2\chi _{h}+4\leq n.$ Taking the reciprocal of this and
combining it with
 (\ref{s-and-chi_h}) yields the result.
\end{proof}

\begin{rem} \textup{The least number of
vanishing cycles for a genus $  2$ Lefschetz fibration  has been
narrowed down to a number that is equal to $7$ or $8$, \cite{Oz} .
Remark \ref{multiple-of-five} gives a minimum value for $n$, which
is $4$, as well as Corollary \ref{proof-2-for-genus-2-bound}. With
that value of $n$ the smallest $s$ can be is $3$ by Corollary
\ref{system-of-equality-inequality}. Therefore the fibration with
$n+s=4+3=7$ vanishing cycles constructed by Xiao in \cite{X}
realizes that minimum number.}
\end{rem}

 From geographical perspective there are three important
regions for genus $2$ Lefschetz fibrations that are distinct in some
ways from one another:

\smallskip

\begin{enumerate}
  \item $2\leq \lambda \leq 4, $
  \item $ 4<\lambda <5, $
  \item $ 5 \leq \lambda < 6$.
\end{enumerate}

\smallskip

In the first region we see most of the known genus $2$ Lefschetz
fibrations that come from topological constructions and mapping
class group considerations. These are the fibrations satisfying
$0\leq\frac{s}{n}\leq\frac{1}{3}$. In particular $\lambda=2$
corresponds to the classical examples that do not contain any
separating vanishing cycles. $\lambda = 4$ corresponds to the
fibrations satisfying $3s=n.$ The well known construction by
Matsumoto has been the only known example satisfying this ratio. The
author of this article has recently given many more examples
satisfying $3s=n.$

 The second region is the loci of fibrations
satisfying $\frac{1}{3}<\frac{s}{n}<\frac{3}{4}$. To the best of our
knowledge there are no known examples of genus $2$ Lefschetz
fibrations in this region coming from topological constructions or
mapping class group considerations. The author of this article has
constructed an example with $\frac{s}{n}=\frac{17}{36}$. All genus
$2$ Lefschetz fibrations in the first two regions satisfy $2n-s>5$
because Proposition \ref{bound-on-slope-is-sharp-prop} requires
$\lambda = 6-\frac{1}{\chi _{h}+1}$ for fibrations satisfying
$2n-s=5$ and $6-\frac{1}{\chi _{h}+1} \geq 5$.

The third region is the region of fibrations satisfying
$\frac{3}{4}\leq \frac{s}{n}<2.$ The fibrations  satisfying the
relation $2n-s=5$ are   in this region. The only known, to the
author,  examples of this sort come from algebro-geometric
constructions and are due to Xiao, \cite{X}. They correspond to
ratios $\frac{s}{n}=\frac{3}{4},\frac{7}{6},\frac{19}{12}$. It's an
open question how high this ratio can be. It would also be
interesting to find a fibration in this region with $2n-s>5$ that is
not a fiber sum of fibrations satisfying $2n-s=5$.

\section{Summary of genus $3$ case}

Almost all of the calculations in the previous section can be
carried out for genus $3$ in much the same manner. We will just list
the results in the sequence they appeared for genus $2$ instead of
redoing all of them.\\

Formula (\ref{slope-formula}) gives
\begin{eqnarray} \label{slope-formula-genus-3} \lambda=4\frac{2n+17s}{3n+8s}=4\frac{2+17r}{3+8r}
\end{eqnarray}
when we substitute $g=3,x=2s$. \\

Proposition \ref{prop-burak}\ \textup{(Corollary 10, \cite{Oz})}
becomes \[c_{1}^{2}\leq \frac{29}{4}\chi
_{h}-\frac{11}{4}=\allowbreak 7.\,25\chi _{h}-2.\,75.\]

Remark \ref{holom-nonneg} becomes $-1\leq \chi _{h}$.\\

Corollary \ref{slope-estimate} gives \quad
\begin{eqnarray} \label{slope-estimate-genus-3}\lambda
\leq \allowbreak \frac{29}{4}-\frac{5}{4}\frac{1}{\chi _{h}+2}.
\end{eqnarray}

Corollary \ref{genus-2-rho} turns out to be \quad $\displaystyle
\rho \left( 3\right) \leq \frac{11}{8}=1.\,375$.

 Corollary \ref{system-of-equality-inequality} takes the form
\begin{eqnarray*}3n+8s&=&28k\\ 11n-8s&\geq&28
\end{eqnarray*}
and solving the system with equalities gives $
s=-\frac{3}{4}+\frac{11}{4}k,n=2+2k,k=\chi _{h}+2\in\mathbb{Z}^+$.
After letting $k=4m+1,m\geq 0,$ we obtain
\begin{eqnarray} \label{ratio-genus-3} \frac{s}{n}=\frac{11m+2}{8m+4},
\end{eqnarray}
which is the genus $3$ version of (\ref{slope-using-m}). Combining
(\ref{slope-formula-genus-3}) and (\ref{slope-estimate-genus-3}) and
using  \[\chi _{h}=\frac{3}{28}n+\frac{2}{7}s-2 \ \ \textup{and} \ \
11n\allowbreak -8s=28\] together we see that the bound
(\ref{slope-estimate-genus-3}) on $\lambda$ would be sharp if there
were fibrations satisfying the equation $11n\allowbreak -8s=28$ but
we do not know any
  example of that. For such fibrations the signature bound
(\ref{burak-sign-bound}) would also be sharp and  realized by genus
$3$ hyperelliptic Lefschetz fibrations satisfying $11n\allowbreak
-8s=28$:
\begin{eqnarray*}11n\allowbreak -8s&=&28\\
\allowbreak n-s-4&=&24-10n+7s\\
&=&\allowbreak 24-10n+7\left( \frac{11}{8}n-\frac{7}{2}\right) \\
&=&\allowbreak -\frac{3}{8}n-\frac{1}{2}\\
&=&\allowbreak -\frac{4}{7}n+\frac{1}{7}\left(
\frac{11}{8}n-\frac{7}{2}\right)\\
&=&-\frac{4}{7}n+\frac{1}{7}s\\
&=&\sigma.
\end{eqnarray*}
In fact, $11n-8s$ is divisible by $28$:
\begin{eqnarray} \label{genus-3-divisibility} 11n-8s=28t,
\end{eqnarray}
where $t=\frac{1}{4}\left( n-s-\sigma \right)\in \mathbb{Z}^{+}$ and
the calculation above is just $t=1$ case (See Remark
\ref{t-general}). Solving a similar system as in Remark
\ref{multiple-of-five} gives
\[\lambda =\allowbreak \frac{29}{4}-\frac{5}{4}\frac{t}{\chi _{h}+2}\leq \frac{29}{4}-\frac{5}{4}\frac{1}{\chi _{h}+2}.\]
\indent  Corollary \ref{slope-bound-in-n-and-s-only} would take the
form

\begin{eqnarray*}\lambda =4\frac{2n+17s}{3n+8s}&\leq& \frac{29s+8}{4s+3}\leq
\frac{1}{4}\frac{87n+232s-140}{3n+8s}\\
&\leq& \frac{1}{4}\frac{29n-68}{n-2}\leq 2\frac{15n+26s-28}{3n+8s}
\leq 2\frac{5n-s-12}{n-2}.
\end{eqnarray*}

\indent All but the fourth inequality above are equivalent to $
0\leq 11n-8s-28  $. The fourth one comes down to $0\leq \left(
3n-16\right) \left( 11n-8s-28\right) $ but $n\geq 8$ for genus $3$
hyperelliptic Lefschetz fibrations. Remark \ref{signature-negative}
would still be  valid for genus $3$ hyperelliptic Lefschetz
fibrations.

Genus $3$ equivalent of the bounds in Corollary
\ref{cor-signature-holom-bound} are \[\displaystyle \sigma \leq
-\frac{3}{4}\chi _{h}-\frac{11}{4} \  \textup{and} \   \sigma \leq
-\frac{3}{19}\chi -\frac{44}{19}.\]

Finally, genus $3$ version of Corollary
\ref{proof-2-for-genus-2-bound} is
\begin{eqnarray*} \frac{s}{n}\leq \frac{11\chi _{h}+19}{8\left( \chi _{h}+3\right) }
\end{eqnarray*}
using $\displaystyle s\leq \frac{1}{4}\left( 11\chi _{h}+19\right),$
which is equivalent to $11n-8s\geq 28$, and $2\chi _{h}+6\leq n$,
(\ref{n-g-bound-on-holom}).

\section{General Case}

\begin{prop}\label{prop-slope-general-formula} For a    genus  $g$ hyperelliptic Lefschetz fibration
the slope is given by  \begin{eqnarray}
\label{slope-general-formula} \lambda =12-\frac{n+s}{\chi
_{h}+g-1}.\end{eqnarray}
\end{prop}
\begin{proof}
 By definition
 \begin{eqnarray*}\lambda &=&\frac{c_{1}^{2}+8\left( g-1\right) }{\chi
 _{h}+g-1}=\frac{12\chi _{h}-\chi +8\left( g-1\right) }{\chi
 _{h}+g-1}\\
 &=&\frac{12\chi _{h}+12g-12-\chi -4\left( g-1\right) }{\chi
 _{h}+g-1}\\
 &=&12+\frac{-\left( n+s-4\left( g-1\right) \right) -4\left( g-1\right) }{\chi
 _{h}+g-1}\\
&=&12-\frac{n+s}{\chi _{h}+g-1}
 \end{eqnarray*}
 \end{proof}

\begin{rem} \textup{ To see that (\ref{slope-general-formula})  agrees with
(\ref{slope-formula}) and (\ref{slope-formula-genus-3})  for genus
$2$ and $3$  simply substitute $\frac{1}{10}\left( n+2s\right) -1$
and $\frac{3}{28}n+\frac{2}{7}s-2$ for $\chi_h$, respectively. The
proof when $s=0$
is straightforward : \\
\[\chi _{h}+g-1=\frac{1}{4}\left( -\frac{g+1}{2g+1}n+n-4\left( g-1\right) \right)
+g-1=\frac{1}{4}\frac{ng}{2g+1}\] and
\[12-\frac{n}{\chi _{h}+g-1}=12-\frac{n}{\frac{1}{4}n\frac{g}{2g+1}}=\allowbreak
4\frac{g-1}{g}.\] }
\end{rem}
\begin{rem} \label{alterbate-general-slope-formulas}\textup{ (\ref{slope-general-formula}) can also be written as
\begin{eqnarray*}\lambda =12-\frac{4}{1+\frac{\sigma }{n+s}}=12-4\frac{n+s}{\sigma +n+s}=8+4\frac{\sigma }{\sigma
+n+s},
\end{eqnarray*}
either by solving the formula given by Theorem
\ref{signature-in-terms-of-lambda} for $\lambda$ or using the
relation
\begin{eqnarray} \label{n+s+sigma-divisible-by-4} \sigma +n+s=4\left(
\chi_h+g-1\right).
\end{eqnarray}
}
\end{rem}
\begin{rem} \label{the-slope-is-less-than-ten-rem} \textup{ The first formula in Remark
\eqref{alterbate-general-slope-formulas} shows how the slope depends
on the (unweighted) "average $\frac{\sigma }{n+s}$ of signature per
vanishing cycle". When $\lambda=10$, this average must be $1$. This
can never happen because the "signature contribution" of each
vanishing cycle is either $-1,$ or $0$, or $+1$ and according to the
handlebody decomposition of Lefschetz fibrations the first handle
attached along the first vanishing cycle, which can be arranged to
be a non-separating one by cyclically permuting, will always result
in a $4-$ manifold with $0$ signature, \cite{Oz}. This is proved in
the following proposition.}
\end{rem}
\begin{prop} \label{slope-estimate-general-prop} For a genus $g$ hyperelliptic Lefschetz fibration we
have
\begin{equation}\label{slope-estimate-general} \lambda \leq 10-\frac{2+s}{\chi _{h}+g-1}.
\end{equation}
\end{prop}

\begin{proof}
First we estimate  $\chi_h$ as
 \begin{eqnarray} \label{definition-of-M} \nonumber \chi _{h}&=&\frac{1}{4}\left( \sigma +\chi \right) =
 \frac{1}{4}\left( -\frac{g+1}{2g+1}n+\sum_{h=1}^{\left[ g/2\right] }\frac{4h\left( g-h\right) s_{h}}{2g+1}-s+n+s-4\left( g-1\right)
 \right)\\ \nonumber
 &\leq& \frac{1}{4}\left( \frac{ng}{2g+1}+\frac{4\frac{g}{2}\left( g-\frac{g}{2}\right) s}{2g+1}-4\left( g-1\right)
 \right)\\
 & =&\frac{1}{4}\frac{\allowbreak ng}{2g+1}+\frac{1}{4}\frac{sg^{2}}{2g+1}-\left( g-1\right) :=M,
 \end{eqnarray} using the fact that $h(g-h)\leq \frac{g}{2}(g-\frac{g}{2})$ and $\sum_{h=1}^{\left[ g/2\right] }s_h=s$. Now,  use this to write $\chi$
 as
 \begin{eqnarray} \label{euler-in-terms-of-K}\nonumber
 \chi &=&n+s-4\left( g-1\right)\\ \nonumber
 &=&\frac{4\left( 2g+1\right) }{g}\left( \allowbreak \frac{1}{4}\frac{ng}{2g+1}+\frac{1}{4}\frac{sg^{2}}
 {2g+1}-\left( g-1\right) \right) +\left( 1-g\right) s+4g-\frac{4}{g}\\
 &=&\frac{4\left( 2g+1\right) }{g}M+\left( 1-g\right)
 s+4g-\frac{4}{g} .
  \end{eqnarray} The estimate
 \begin{eqnarray*}
  \sigma &\leq& n-s-4=n+s-4\left( g-1\right) -2s+4\left( g-2\right) =\chi -2s+4\left(
  g-2\right),
  \end{eqnarray*}
  (\ref{burak-sign-bound}), can be used to write
 \begin{eqnarray} \label{holom-bdd-by-euler}
 \hspace{-.35in} \chi _{h}&=&\frac{1}{4}\left( \sigma +\chi \right) \leq \frac{1}{4}\left( \chi -2s+4\left( g-2\right) +\chi \right) =\frac{1}{2}\chi
  -\frac{1}{2}s+g-2
 \end{eqnarray} and using (\ref{euler-in-terms-of-K}) we obtain
\begin{eqnarray*}\chi _{h}&\leq &\frac{1}{2}\left( \frac{4\left( 2g+1\right) }{g}M+\left( 1-g\right) s+4g-\frac{4}{g}\right)
-\frac{1}{2}s+g-2\\
 &=&\allowbreak 2\frac{2g+1}{g}M-\frac{1}{2}sg+3g-2-\frac{2}{g}.
 \end{eqnarray*}
 We will solve this for $sg$
  \begin{eqnarray*}sg\leq \allowbreak 4\frac{2g+1}{g}M-2\chi _{h}+6g-4-\frac{4}{g}\end{eqnarray*}
 and use it in estimating
 \begin{eqnarray*}
 c_{1}^{2}&=&12\chi _{h}-\chi =12\chi _{h}-\left( \frac{4\left( 2g+1\right) }{g}M+\left( 1-g\right)
 s+4g-\frac{4}{g}\right)\\
 &=&12\chi _{h}-4\frac{2g+1}{g}M+\left( g-1\right) s-4g+\frac{4}{g}\\
 &\leq & 12\chi _{h}-4\frac{2g+1}{g}M+\allowbreak 4\frac{2g+1}{g}M-2\chi
 _{h}+6g-4-\frac{4}{g}-s-4g+\frac{4}{g}\\
 &=&10\chi _{h}+2g-4-s.
 \end{eqnarray*}
 Now,
 \begin{eqnarray*}\lambda =\frac{c_{1}^{2}+8\left( g-1\right) }{\chi _{h}+g-1}\leq \frac{10\chi _{h}+2g-4-s+8\left( g-1\right) }{\chi _{h}+g-1}=\frac{10\chi _{h}+10g-10-2-s}{\chi _{h}+g-1}
 \end{eqnarray*}
 and we have
 \begin{eqnarray*}\lambda \leq 10-\frac{2+s}{\chi _{h}+g-1}.
 \end{eqnarray*}
\end{proof}
\begin{cor} \label{the-slope-is-less-than-ten-cor} The slope $\lambda$ of an hyperelliptic genus $g$
Lefschetz fibration satisfies $\lambda \leq 10$.
\end{cor}
\begin{rem} \label{generalizations-rem} \textup{ Proposition \ref{prop-slope-general-formula} is true in
general, i.e., the assumption that the Lefshcetz fibration  is
 hyperelliptic is not necessary. Therefore the formulas in Remark
\ref{alterbate-general-slope-formulas} are also true in general and
using Remark \ref{the-slope-is-less-than-ten-rem} we can say that
Corollary \ref{the-slope-is-less-than-ten-cor} extends to
non-hyperelliptic Lefschetz fibrations as well. Because of Remark
\ref{alterbate-general-slope-formulas} we also conclude that Theorem
\ref{signature-in-terms-of-lambda} extends to non-hyperelliptic
fibrations.  }
\end{rem}
\begin{rem} \textup{ \label{min-n-is-4} One can show that
\begin{eqnarray} \label{n-g-bound-on-holom}
 \chi _{h}\leq\frac{n}{2}-g, \ \ \ \textup{i.e.,} \ \ \ 2\chi _{h}+2g\leq n,
 \end{eqnarray}
for hyperelliptic genus $g$ Lefschetz fibrations using
(\ref{holom-bdd-by-euler}):
\begin{eqnarray*}
  \chi _{h}\leq \frac{1}{2}\chi
  -\frac{1}{2}s+g-2=\frac{1}{2}\left( n+s-4\left( g-1\right) \right)
  -\frac{1}{2}s+g-2=\frac{1}{2}n-g.
\end{eqnarray*} }
\end{rem}
\begin{cor}  Let $X
\rightarrow S^{2}$ be a simply connected genus $g\geq 2$
hyperelliptic Lefschetz fibration with $b_{2}^{+}\geq 1.$ Then the
minimum number of non-separating vanishing cycles is $2g+2$. If
furthermore $b_{2}^{+}> 1$  then this minimum becomes $2g+4$.
\end{cor}
\begin{proof} By definition of  $\chi _{h}$ we have
 \[\chi _{h}=\frac{1}{4}\left( \sigma +\chi \right) =\frac{1}{4}\left(
b_{2}^{+}-b_{2}^{-}+2-2b_{1}+b_{2}^{+}+b_{2}^{-}\right)
=\frac{1}{2}\left( b_{2}^{+}+1-b_{1}\right). \] Using
\eqref{n-g-bound-on-holom} and the assumption $b_1=0$ we get
\[\frac{1}{2}\left( b_{2}^{+}+1\right) \leq \frac{1}{2}n-g.\]
Solving this inequality for $n$ after using  $b_{2}^{+}\geq 1$
yields $2g+2\leq n.$ Clearly   $2g+4\leq n$ when $b_{2}^{+}> 1$
because $b_{2}^{+}$ must be odd.
\end{proof}

\begin{prop} \label{slope-double-estimate-prop} The slope of an hyperelliptic genus $g$
Lefschetz fibration satisfies
\begin{eqnarray} \label{slope-double-estimate}
 4\frac{g-1}{g}+\frac{4s}{g}\cdot\frac{\left( 2g+1\right) \left( 3g-4\right) }{ng+4s\left( g-1\right) }\leq \lambda \leq 10-2\frac{2+s}{n-2}.
\end{eqnarray}
\end{prop}
\begin{proof} The signature satisfies the bound
\begin{eqnarray*}  \sigma =-\frac{g+1}{2g+1}n+\frac{4x}{2g+1}-s\geq -\frac{g+1}{2g+1}n+\frac{4s(g-1)}{2g+1}-s
=-\frac{g+1}{2g+1}n+\frac{2g-5}{2g+1}s
\end{eqnarray*}
because $s(g-1)\leq x$ by definition of $x$ and $s$. Now, using
Theorem \ref{signature-in-terms-of-lambda} we can write
\[-\frac{g+1}{2g+1}n+\frac{2g-5}{2g+1}s\leq -\frac{8-\lambda }{12-\lambda }\left( n+s\right) \]
and solving this for $\lambda$ gives the first inequality. To prove
the second inequality we begin with the fact that $\chi _{h}+g-1>0,$
as we mentioned in the proof of Corollary \ref{min-n-is-4}. Using
this and (\ref{n-g-bound-on-holom}) we can write
\[-\frac{1}{\chi _{h}+g-1}\leq \allowbreak \frac{-2}{n-2}.\]
Now, adding $10$ to both sides after multiplying by $2+s$  proves
the second inequality thanks to Proposition
\ref{slope-estimate-general-prop}.
\end{proof}
\begin{rem} \textup{ We wrote (\ref{slope-double-estimate}) in that particular form   instead of simplifying it
in order to emphasize the fact that it is another proof for Theorem
\ref{main-thm} and that $4-\frac{4}{g}\leq \lambda \leq 10$ for
hyperelliptic Lefschetz fibrations. The lower bound in
(\ref{slope-double-estimate}) gives (\ref{slope-in-n-and-s}) when we
set $g=2$ and it gives the genus $3$ version of
(\ref{slope-in-n-and-s}) when $g$ is set equal to $3$. The reason
this estimate is sharp for low genus is the fact that there is only
one type of separating vanishing cycle for low genus and due to that
reason the estimate $s(g-1)\leq x$ becomes equality for genus
$g=2,3.$ }
\end{rem}
\begin{prop} \label{n-even} Let  $X \rightarrow S^{2}$ be a genus $g$  hyperelliptic Lefschetz fibration
with $n$ non-separating vanishing cycles. Then \\

\indent
  \begin{tabular}{ll}
   $\bullet\ \, n$ \ \textup{is divisible by} $4$, & \hbox{\textup{if} $g$ \textup{is odd};} \\
   $\bullet \ \, n$ \ \textup{is even}, & \hbox{\textup{if} $g\equiv 2 $ (\textup{mod} $4$)  .}
  \end{tabular}
\end{prop}

\begin{proof} $\sigma +s+n$ is divisible by $4$ by (\ref{n+s+sigma-divisible-by-4}).
 Write the signature
\[\sigma =-\frac{g+1}{2g+1}n+\frac{4x}{2g+1}-s,\]
where $x=\sum_{h=1}^{\left[ g/2\right] }h\left(g-h\right) s_{h}$, as
\begin{eqnarray} \label{gn-divisible-by-4} \left( 2g+1\right) \left( \sigma +s\right) +\left( g+1\right) n=4x.
\end{eqnarray}
Equivalently,
\begin{eqnarray*}  \left( 2g+1\right) \left( \sigma +s+n\right) -gn=4x,
\end{eqnarray*}
which shows that $gn$ is divisible by $4$ and the proof follows from
that.
\end{proof}
Divisibility of $n$ by $4$ when $g$ is odd also follows from
Proposition 4.10 of \cite{E}.

\begin{rem} \label{t-general}\textup{ If $g$ is not divisible by $4$ then $n$ is even by Proposition \ref{n-even}.
In that case we   conclude from (\ref{gn-divisible-by-4}) that
$s+\sigma$ is also even. We use this and the fact that   $\sigma
+s+n$ is divisible by $4$ to prove that $n-s-\sigma $ is divisible
by $4$ as well when  $g$ is not divisible by $4$:
\begin{eqnarray*} n-s-\sigma &=&\sigma +s+n-2\left( s+\sigma \right).
\end{eqnarray*}
Then
\begin{eqnarray}\label{general-divisibility} \frac{1}{4}\left( n-s-\sigma \right) =
\allowbreak \frac{1}{4}\left( n-s-\left(
-\frac{g+1}{2g+1}n+\frac{4x}{2g+1}-s\right) \right) =
\frac{1}{4}\frac{\left( 3g+2\right) n-4x}{2g+1}\in \mathbb{Z}^+.
\end{eqnarray}
When $g=2$, (\ref{general-divisibility}) becomes $2n-s=5\left(
n+\sigma \right) $ as we found in Remark \ref{multiple-of-five}.
When $g=3$ then (\ref{general-divisibility}) is the same as
(\ref{genus-3-divisibility}). The integer
(\ref{general-divisibility}) is positive because of
(\ref{burak-sign-bound}).}
\end{rem}
\begin{proof}[Proof of Theorem \ref{ratio-estimate}]
Using the bound  (\ref{burak-sign-bound}) we get $\displaystyle
\frac{1}{4}\left( n-s-\sigma \right)\geq 1$. Then
(\ref{general-divisibility}) gives
\begin{eqnarray*} 1\leq \allowbreak \frac{1}{4}\frac{\left( 3g+2\right) n-4x}{2g+1},
\end{eqnarray*}
and hence
\begin{eqnarray*} \allowbreak x\leq \frac{1}{4}n\left( 3g+2\right) -\left( 2g+1\right).
\end{eqnarray*}
 Using the estimate  $(g-1)s\leq x$ one more time, we have
 \begin{eqnarray*}\left( g-1\right) s\leq \frac{1}{4}n\left( 3g+2\right) -\left( 2g+1\right).
\end{eqnarray*}
Dividing through by $n(g-1)$ gives
\begin{eqnarray*}r=\frac{s}{n}\leq \frac{3g+2}{4\left( g-1\right) }-\frac{2g+1}{n\left( g-1\right) }.
\end{eqnarray*}
Since $s$ and $n$ are arbitrary, we conclude
\[\rho \left( g\right) \leq \frac{3g+2}{4\left( g-1\right) }.\]
\end{proof}
\begin{cor} For an hyperelliptic Lefschetz fibration of genus $g\geq 6$ we have $s\leq n$.
\end{cor}
\begin{rem} \label{final-remarks-and-alternate-proofs} \textup{ One can   prove Theorem \ref{ratio-estimate} by solving
\begin{eqnarray*}
 4\frac{g-1}{g}+\frac{4s}{g}\cdot\frac{\left( 2g+1\right) \left( 3g-4\right) }{ng+4s\left( g-1\right) } \leq
 10-2\frac{2+s}{n-2}
\end{eqnarray*}
 for $\frac{s}{n}$ as well, \eqref{slope-double-estimate}.
Also, solving
\begin{eqnarray*}\lambda =12-4\frac{n+s}{n+s+\sigma }\leq
10-2\frac{2+s}{n-2}
\end{eqnarray*}
for $\sigma$ results in \eqref{burak-sign-bound}, which is another
proof  for Proposition \ref{slope-estimate-general-prop}. Finally,
solving
\begin{eqnarray*}4\frac{g-1}{g}\leq \lambda=12-\frac{4}{1+\frac{\sigma }{n+s}}
\end{eqnarray*} for $\frac{\sigma }{n+s}$ gives
\begin{eqnarray*}\frac{\sigma }{n+s}\geq -\frac{g+1}{2g+1},
\end{eqnarray*} which shows that "the average signature per vanishing cycle"   is
at least $-\frac{g+1}{2g+1}$ for Lefschetz fibrations satisfying
$\lambda \geq 4-4/g$ and it is greater than that whenever $s>0$ by
virtue of Theorem \ref{main-thm}.\\
 Based on this observation we
conclude the following bound on $ \rho \left( g\right)$ in general,
without assuming hyperellipticity:}
\end{rem}

\begin{cor} For a Lefschetz fibration  of genus $g\geq 2$ we have
\[\rho \left( g\right) < 3+\frac{2}{ g }\, .\]
\end{cor}
\begin{proof} By Corollary 7 in \cite{Oz} we have $\sigma \leq n-s$.
Combining that with Remark \ref{final-remarks-and-alternate-proofs}
we conclude
\begin{eqnarray*} -\frac{g+1}{2g+1}<
\frac{n-s}{n+s}=\frac{1-r}{1+r}.
\end{eqnarray*}
 The result follows once we solve this inequality for
$r$.
\end{proof}

\emph{Acknowledgements.}\label{ackref} Many thanks to Hur\c{s}it
\"{O}nsiper for helpful and encouraging conversations and for
referring the author to the article written by Xiao. The author is
also grateful to Hisaaki Endo for his insightful comments.

   \noindent Yusuf Z. G\"{u}rta\c{s}\\
  Mathematics Department\\
   DePauw University\\
602 S. College Avenue\\
Greencastle, IN 46135\\
U.S.A. \\
  yusufgurtas@depauw.edu


%
\end{document}